\documentclass{proc-l}
\usepackage{amssymb} % Only used it at \nmid

\newcommand{\F}{{\mathbb F}}
\newcommand{\Z}{{\mathbb Z}}

\newcommand{\cS}{{\mathcal S}}

\newcommand{\cH}{{\mathcal H}}
\def\P{\mathbb P}

\DeclareMathOperator{\Aut}{Aut}
\DeclareMathOperator{\End}{End}
\DeclareMathOperator{\Gal}{Gal}

\newtheorem{theorem}{Theorem}
\newtheorem{lemma}[theorem]{Lemma}
\newtheorem{corollary}[theorem]{Corollary}
\newtheorem{question}[theorem]{Question}

\theoremstyle{remark}
\newtheorem{case}{Case}
\newtheorem{acks}{Acknowledgments}

\begin{document}

\title[Hyperelliptic curves without extra automorphisms]
{Hyperelliptic curves over $\F_2$ of every $2$-rank without extra
automorphisms}
\author{Hui June Zhu}
\address{Department of Maths and Stats, McMaster University,
Hamilton, ON L8S 4K1, CANADA}
\email{zhu@cal.berkeley.edu}
\date{March 28th, 2005}
\commby{David E. Rohrlich}

%%%%%%%%%%%%%%%%%%%%%%%%%%%%%%%%%%%%%%%%%%%%%%%%%%%%%%%%%%%%%%%%%%%%%%%%%%%%%
%%                                                                         %%
%% (2000 Mathematics Subject Classification)                               %%
%%                                                                         %%
%% 11            Number theory                                             %%
%%    11G           Arithmetic algebraic geometry (Diophantine geometry)   %%
%%        11G10        Abelian varieties of dimension $\gtr  1$            %%
%% 14            Algebraic Geometry                                        %%
%%    14G           Arithmetic problems. Diophantine geometry              %%
%%        14G15        Finite ground fields                                %%
%%    14K           Abelian varieties and schemes                          %%
%%        14K15        Arithmetic ground fields                            %%
%%                                                                         %%
%%%%%%%%%%%%%%%%%%%%%%%%%%%%%%%%%%%%%%%%%%%%%%%%%%%%%%%%%%%%%%%%%%%%%%%%%%%%%
\begin{abstract}
We prove that for any pair of integers $0\leq r\leq g$ such that
$g\geq 3$ or $r>0$, there exists a (hyper)elliptic curve $C$ over $\F_2$
of genus $g$ and $2$-rank $r$ whose automorphism group consists of
only identity and the (hyper)elliptic involution.
As an application, we prove the existence of principally
polarized abelian varieties $(A,\lambda)$ over $\F_2$ of dimension $g$
and $2$-rank $r$ such that $\Aut(A,\lambda)=\{\pm 1\}$.
\end{abstract}

\maketitle

\section{Introduction}

In this paper curves are smooth, projective, and geometrically
integral algebraic varieties of dimension one defined over fields.
Let $k$ be a field and $\overline{k}$ its algebraic closure. If
$C$ is a curve over $k$, let $\Aut{C}$ denote the group of
automorphisms of $C$ defined over $\overline{k}$. Let $J(C)$
denote the Jacobian of $C$. Let $\End{J(C)}$ denote the
endomorphism ring of $J(C)$ over $\overline{k}$. Let $\F_p$ be a
finite field of $p$ elements for some prime $p$. Let
$\overline\F_p$ be its algebraic closure.

A supersingular curve $C$ over $\F_p$ is a curve whose Jacobian is
isogenous over $\overline\F_p$ to a product of supersingular
elliptic curves. Hence a supersingular curve $C$ is a cover of
these supersingular elliptic curves. It has $p$-rank $0$ but the
converse is not true for $g\geq 3$. Supersingular curves are
intimately connected to curves with large automorphism groups. For
instance, in the seminal paper~\cite{GeerVlugt}, the authors
constructed supersingular curves over finite field of
characteristic $2$ by taking quotients of some families of
($2$-rank $0$) curves over $\F_2$ with large automorphism groups.
It is well-known that curves over fields of positive
characteristic achieving maximal automorphism groups are all
supersingular curves~\cite{Stichtenoth}. Is it a myth or truth
that a curve over $\F_p$ of lower $p$-rank has larger automorphism
groups in general?

In the moduli space of curves, the subset corresponding to the curves
with trivial automorphism group is open
(see~\cite[Introduction]{Poonen1} or~\cite[Remark
10.6.24]{Katz-sarnak}). In a recent paper this fact was proved
constructively~\cite{Poonen1} (see also~\cite{Poonen2}\cite{Poonen3}).
It is desirable to understand how this subset stratifies
by the $p$-rank of the curves.

\begin{question}
Let $p$ be a prime number.  Given integers $g\geq 3$ and $0\leq r\leq
g$, is there a curve $C$ over $\F_p$ of genus $g$ and $p$-rank $r$ such
that $\Aut C=\{1\}$?
\end{question}

There is not any constructive way to obtain curves
over $\F_p$ of prescribed genus and $p$-rank, so we do not know the
answer to this question.

On the other hand, for every prime $p$ and positive integer $g$,
Poonen~\cite{Poonen2} has constructed (hyper)elliptic curves $C$ over
$\F_p$ of genus $g$ with $\Aut C=\{1,\iota\}$, where $\iota$ is the
unique (hyper)elliptic involution of $C$. Automorphisms other than
these two are referred as {\em extra} automorphisms.

If $g=1$, it is well-known that for every prime $p$ a supersingular
elliptic curve (i.e., with zero $p$-rank) over $\F_p$ has extra
automorphisms, while there exist ordinary elliptic curves (i.e., with
non-zero $p$-rank) over $\F_p$ with $\Aut C=\{1,\iota\}$.
(See~\cite[Chapter~III]{Silverman}.)

Every curve over $\F_2$ of genus $2$ and $2$-rank $0$ can be written
in the form $y^2+y = x(x^4+a_1x^2+a_0x)$ for $a_0,a_1\in\F_2$, hence
has extra automorphisms. It is easy to check this fact by hand. In
fact, every curve of the form $y^2+y = x(\sum_{i=0}^{n}a_ix^{2^i})$
for some integer $n$ and $a_i\in\F_2$ has extra automorphisms
(see~\cite{GeerVlugt}).

\begin{question}
Let $p$ be a prime number. Given integers $g\geq 1$ and $0\leq r\leq
g$, is there a (hyper)elliptic curve $C$ over $\F_p$ of genus $g$ and
$p$-rank $r$ without extra automorphism?
\end{question}

The present paper gives a complete answer to this question for the
case $p=2$. We hope this provides evidence for more general
theorem or conjecture in the future.

\begin{theorem}
For any integers $0\leq r\leq g$ such that $g\geq 3$ or $r>0$, there
exists a (hyper)elliptic curve $C$ over $\F_2$ of genus $g$ and
$2$-rank $r$ such that $\Aut(C)=\{1,\iota\}$ where $\iota$ is the
unique (hyper)elliptic involution of $C$.
\end{theorem}

The proof of the theorem, divided in two parts, is presented in
the next two sections. This theorem has the following application.
For an abelian variety $A$ with polarization $\lambda$ defined
over $\F_2$ let $\Aut(A,\lambda)$ denote the group of
automorphisms of $A$ over $\overline\F_2$ respecting the
polarization. The corollary below follows immediately from the
theorem by applying the Torelli's theorem~\cite[Theorem
12.1]{Milne}.  Detailed discussion upon related results can be
found in the Introduction of~\cite{Poonen2}.

\begin{corollary}
For any integers $0\leq r\leq g$ such that $g\geq 3$ or $r>0$, there
exists a $g$-dimensional principally polarized abelian variety
$(A,\lambda)$ over $\F_2$ of $2$-rank $r$ such that
$\Aut(A,\lambda)=\{\pm 1\}$.
\end{corollary}

Finally we remark that these above two questions will be resolved
if we know what algebras can become $\End{J(C)}$ for curves $C$ over
$\F_p$ of prescribed genus (see~\cite[Question (8.6)]{Oort}) and
$p$-rank. By \cite{Mori} (see also \cite{Zarhin}), one knows that
for every $g\geq 1$ there exists a hyperelliptic curve $C$ of any
genus $g\geq 1$ with $\End{J(C)}=\Z$. However, this does not hold
for curves over finite fields, in which case we have $\End{J(C)}$
strictly contains $\Z$.

\section{Construction for $r>0$}

Suppose $g\geq 2$ and $r\leq g$ are two positive integers. Let
$q(x)$ be a polynomial in $\F_2[x]$ of degree $<2g+1-r$
(resp. $=2g+1-r$) with $r$ (resp. $r+1$) distinct roots and, let
$f(x)$ be a polynomial in $\F_2[x]$ of degree $2g+1-r$ (resp. $\leq
2g+1-r$), such that $f(x)$ and $q(x)$ has no common roots.
Let $C$ be the hyperelliptic curve over $\F_2$ defined by the affine equation
\begin{equation}\label{C1}
C: y^2+y=\frac{f(x)}{q(x)}.
\end{equation}

Then the curve $C$ over $\F_2$ is of genus $g$ by the Riemann-Hurwitz
formula and of $2$-rank $r$ by the Deuring-Shafarevich formula in
(\ref{C3}), which we shall explain immediately (see details
in~\cite{Madan} or~\cite{Nakajima}).  Let $k$ be an algebraically closed
field of characteristic $p$. Let $\pi:X\rightarrow Y$ be a finite
Galois covering of curves over $k$ whose Galois group $G$ is a
$p$-group. Let $r_X$ and $r_Y$ denote the $p$-ranks of $X$ and $Y$,
respectively.  Let $Q_1,\cdots, Q_n$ be the set of ramification points
on $Y$ with respect to $\pi$. For each point $Q_i$ let $p^{e_i}$ (here
$e_i\geq 1$) be its ramification index. Then
\begin{equation}\label{C3}
r_X-1 = \# G\cdot (r_Y-1+\sum_{i=1}^{n}(1-p^{-e_i})).
\end{equation}

Let $D$ be the ramification divisor of the canonical double cover
$C\rightarrow \P^1$. Write $q:=\prod_{i=1}^{r}(x-\alpha_i)^{b_i}$
(resp. $q:=\prod_{i=1}^{r+1}(x-\alpha_i)^{b_i}$) for distinct
$\alpha_i\in \overline\F_2$ and $b_i\in\Z_{>0}$. The set $\cS$ of
ramification points consists of those points $P_{\alpha_i}$
corresponding to the zeroes of $q$ and possibly the point
$P_\infty$ at infinity. We have
\begin{equation}\label{C2}
D =
\left\{
\begin{array}{ll}
(2g+2-\deg(q)-r)P_\infty
   +\sum_{i=1}^{r}(b_i+1)P_{\alpha_i} \quad \mbox{and respectively}\\
   \sum_{i=1}^{r+1}(b_i+1)P_{\alpha_i}
\end{array}\right.
\end{equation}

Every automorphism of $C$ gives rise to an automorphism of $\P^1$
preserving $D$ under the canonical double cover $C\rightarrow \P^1$.
To construct curves $C$ without extra automorphisms, it suffices to
find monic polynomials $f$ and $q$ in $\F_2[x]$ such that every
automorphism of $\P^1$ preserving $D$ is the identity map on $\P^1$.

Our construction below follows the following idea: for every pair of
integers $0<r\leq g$, we shall construct polynomials $q$ such that
$q$ has $r$ (or $r+1$ resp.) distinct roots and
of degree $<2g+1-r$ (or $2g+1-r$, resp.) in $\F_2[x]$. We always let
$f$ be any polynomial in $\F_2[x]$ of degree $2g+1-r$ (or $\leq
2g+1-r$, resp.) which has no common roots with $q$. We remark that
we shall use the construction that $q$ has $r$ distinct roots
except in {\em Case 5} and {\em Case 6}.

In the construction below we use the notation $f_n$ for a $n$-th
degree irreducible polynomial in $\F_2[x]$. It is a basic fact in
algebra that $f_n$ exists for every positive integer $n$
(see~\cite[Chapter V]{Lang}). For example, $f_2= x^2+x+1$ and
$f_3=x^3+x^2+1$ or $x^3+x+1$. For any $f_3$ of our choice, we
denote by $\beta_1,\beta_2,\beta_3$ its roots in $\overline\F_2$
in an order such that $\beta_1^2 = \beta_2$.

\begin{case}
Suppose $r\geq 8$:
\end{case}

Let $q = f_3f_{r-3}$ if $3\nmid r$ and $q=xf_3f_{r-4}$ otherwise.

Let $\sigma$ be an automorphism of $\P^1$ which acts as a
$3$-cycle on the three roots of $f_3$ in $\overline\F_2$. Since
$3$ points determines an automorphism of $\P^1$, $\sigma$ is
defined over the field $k$ generated by roots of $f_3$ over
$\F_2$. Hence, $\sigma(P_\infty)$ corresponds to a point in $k$.

Let $\F$ denote the composition of all finite extensions of $\F_2$ of
degrees coprime to $3$. There are exactly $r-2$ distinct $\F$-rational
points in the set of ramification points $\cS$. Suppose $\lambda$ is a
non-trivial automorphism of $\P^1$ preserving $D$.  Then $\lambda$
must map at least $(r-2)-3\geq 3$ of these $\F$-rational points to
other $\F$-rational points of $\cS$. But $\lambda$ is determined by
its values at $3$ points, so $\lambda$ must be defined over $\F$. In
particular, $\lambda$ preserves the set of $3$ non-$\F$-rational
points of $\cS$, the roots of $f_3$.  If $\lambda$ fixes any one of
them, as they are Galois conjugates over $\F$, then $\lambda$ would
fix them all, hence $\lambda$ would be trivial. So $\lambda$ acts as a
$3$-cycle, and after replacing $\lambda$ by $\lambda^{-1}$ if
necessary, we may assume $\lambda=\sigma$.  Since $\lambda$ permutes
the roots of $f_3$, it fixes its coefficients, hence $\lambda$ fixes
$0$ and $1$. So $\lambda(P_\infty)\neq P_\infty$ and
$\lambda(P_\infty)\neq P_1$.  But $D$ is preserved, so $\lambda$ maps
$P_\infty$ to a root of $f_{r-3}$ (or $f_{r-4}$), which lies in $\F$
and does not lie in $k$. This contradicts our assumption above about
$\sigma$.

\begin{case}
Suppose $r=1$ and $g\geq 2$, or $r=2$ and $g\geq 4$:
\end{case}

For $r=1$ and $g\geq 2$, let $q=x$.

Then the ramification divisor is $D=2gP_\infty + 2P_0$. Since
$g\geq 2$ every automorphism of $\P^1$ preserving $D$ fixes
$\infty$ and $0$, hence it is of the form $x\mapsto cx$ for some
non-zero $c\in \overline\F_2$. A simple computation shows that
$c=1$. This resembles Case I in Section 2 of~\cite{Poonen2}.

For $r=2$ and $g\geq 4$, let $q = x^2(x+1)$.

Then $D=(2g-3)P_\infty+3P_0+2P_1.$
Every automorphism of $\P^1$ preserving $D$ has three points
$\infty, 0$ and $1$ all fixed hence is identity.

\begin{case}
Suppose $r=3$ and $g\geq 4$:
\end{case}

Let $q = f_3$.

Then $D = (2g-4)P_\infty + 2(P_{\beta_1}+P_{\beta_2}+P_{\beta_3})$.
Let $\lambda$ be a non-trivial automorphism of $\P^1$ that preserves
$D$.  By assumption $2g-4>2$, so $\lambda$ fixes $P_\infty$ and
$\lambda$ permutes the roots of $f_3$. Thus $\lambda$ fixes $0$ and
$1$. But then it fixes all three points $0,1$ and $\infty$, it must be
identity. This leads to a contradiction.

These following three cases follow the same scheme, so we shall
elaborate on Case 4 and only sketch the rest two cases.

\begin{case}
Suppose    $r=4$ and $g\geq 5$,
    or $r\geq 4$ and $g\geq r+3$:
\end{case}

For $r=4$ and $g\geq 5$, let $q=x^2f_3$.

Then the ramification divisor is $D=(2g-7)P_\infty
+3P_0+2(P_{\beta_1}+P_{\beta_2}+P_{\beta_3}).$ Let $\lambda$ be a
non-trivial automorphism of $\P^1$ which preserves $D$, then
$\lambda$ permutes the roots of $f_3$ hence it fixes $0$ and $1$.
If fixes $P_\infty$ and $P_0$ then it is identity. If $\lambda$
swaps $P_\infty$ and $P_0$, and it is of the form
$\lambda(\alpha)=c/\alpha$ for some non-zero $c\in \overline\F_2$.
It can be checked quickly that this map can not preserve the roots
of $f_3$.

For $r\geq 4$ and $g\geq r+3$, let $q=x^3(x+1)^2f_{r-2}$.

Then $D = (2g-2r-1)P_\infty + 4P_0 + 3P_1 + 2\sum_{(f_{r-2})_0}P.$
Since $2g-2r-1 \geq 5$ and $r-2\geq 2$, every automorphism of $\P^1$
preserving $D$ has three points $\infty, 0$ and $1$ all fixed hence is
identity.

\begin{case}
Suppose $r=5$ and $g\geq 5$:
\end{case}

Let $q=f_3(x+1)^{2g-9}(x^2+x+1)$. Let $\alpha_1, \alpha_2$ be
roots of $x^2+x+1$ in $\overline\F_2$.

Then the ramification divisor is
$$D=2(P_{\beta_1}+P_{\beta_2}+P_{\beta_3})+(2g-8)P_1+
    2(P_{\alpha_1}+P_{\alpha_2}).$$
Label the roots of $\beta_1,\beta_2,\beta_3,1,\alpha_1,\alpha_2$ by
$1,2,3,4,5,6$, respectively, such that the absolute Frobenius acts on
$\cS$ as the permutation $\sigma = (123)(56)$. Let $H$ be the subgroup
of the automorphism of $\P^1$ preserving $D$, which we may view as a
faithful subgroup of $S_6$, since automorphisms are determined already
by $3$ values. Any automorphism of $\P^1$ which fixes $\alpha_1$ and
$\alpha_2$ has to fix $1$ so $\beta_3$ can not be mapped to $1$.
Therefore, $(12)(34)\not\in H$.
The
group theoretical lemma below, due to Poonen
(see~\cite[Lemma~3]{Poonen2}), indicates that $H$ is trivial.

\begin{lemma}\label{L:poonen}
Suppose $H$ is a subgroup of $S_6$ such that
\begin{enumerate}
\item Each non-trivial element of $H$ has at most $2$ fixed
points; \item $\sigma H \sigma^{-1}\subset H$ for every $\sigma\in
\Gal(\overline\F_2/\F_2)$; \item The permutation $(12)(34)$ is not
in $H$,
\end{enumerate}
Then $H=\{1\}$.
\end{lemma}

\begin{case}
Suppose $r=6$ and $g\geq 7$:
\end{case}

Let $q=f_3(x+1)(x^2+x+1)x^{2g-11}$.

Then the ramification divisor is
$D=2(P_{\beta_1}+P_{\beta_2}+P_{\beta_3})+2P_1
+2(P_{\alpha_1}+P_{\alpha_2})+(2g-10)P_0$.  Note that every
automorphism of $\P^1$ preserving $D$ fixes $P_0$. Then we apply the
same argument as in {\em Case 5}.

\begin{case}
Suppose $r=7$ and $g\geq 8$:
\end{case}

Let $q = f_3(x+1)(x^2+x+1)x^2$.

Then the ramification divisor is
$$D= 2(P_{\beta_1}+P_{\beta_2}+P_{\beta_3}) + 2P_1
+2(P_{\alpha_1}+P_{\alpha_2}) + 3P_0 + (2g-13)P_\infty.$$ Let
$\lambda$ be a non-trivial automorphism of $\P^1$ preserving $D$.
If $\lambda$ fixes $P_0$ and $P_\infty$ then we use the same
argument as in {\em Case 5}.  This is the case when $g\geq 9$.  It
remains to prove the case $g=8$ and $\lambda$ swaps $P_0$ and
$P_\infty$. Then $\lambda(\alpha)=c/\alpha$ for some non-zero
$c\in \overline\F_2$.  If $\lambda$ fixes $P_1$ then it is defined
over $\F_2$ hence it permutes the roots of $f_3$ and fixes $P_0$,
contradiction.  If $\lambda$ swaps $P_1$ with one root of $f_3$,
then it preserves the roots of $f_3$ also. If $\lambda$ swaps
$P_1$ with a root of $f_2$ then it permutes the roots of $f_2$. So
it has to fixes $P_1$, which is absurd.

\begin{case}
Remaining cases:
\end{case}
For $g=r=4,6$, let $C: y^2+y= x+ \frac{1}{x(x^{r-1}+1)}$.

For $g=r=3,5,7$, let $C: y^2+y= x+ \frac{1}{x^r+1}$.

For $g=2,3$ and $r=2$, let $C: y^2+y = x+\frac{1}{x^2+x+1}$ and $C:
y^2+y = x^3+\frac{1}{x^2+x+1},$ respectively.

It is an elementary computation to show that these curves
have no extra automorphisms.

\section{Construction for $r=0$}

We still assume $g\geq 2$.  In this section let $C$ be a
hyperelliptic curve defined by the affine equation
\begin{equation}\label{C4}
C: y^2+y = f(x)
\end{equation}
where $f(x)$ is a polynomial in $\F_2[x]$ of degree $2g+1$.  This is
the same as letting $q=1$ in~(\ref{C1}).  So $C$ is of genus $g$
and $2$-rank $0$. We remark that every curve in (\ref{C4})
is isomorphic to a curve with only odd-degree terms in $f(x)$ because
the base field is $\F_2$.

Any automorphism of $C$ is of the form $x\mapsto ax+b$ and
$y\mapsto cy+h(x)$ for some $a,b,c\in \overline\F_2$ and some
polynomial $h(x)$ in $\overline\F_2[x]$ of $\deg(h)\leq g$.  Let
$\cH$ be the set of polynomials $p(x)^2 + p(x)$ for all polynomial
$p(x)$ in $\overline\F_2[x]$ of degree $\leq g$. It is easy to
show that it is a $\overline\F_2$-vector space of dimension $g+1$.
It follows that $c=a^{2g+1}=1$ and $f(ax+b)+f(x)=h(x)^2+h(x)$.
That is
\begin{equation}\label{C5}
a^{2g+1}=1\quad \mbox{and}\quad f(ax+b)+f(x)\in \cH.
\end{equation}

\begin{lemma}\label{L4}
Let $g=4$ or $g\geq 7$.
Let $p(x)$ be a polynomial in $\F_2[x]$ of degree $\leq
2g-6$. The hyperelliptic curve
$C$ defined by the affine equation
$$
C: y^2 + y = f(x) : = x^{2g+1} + x^{2g-1} + x^{2g-3} + p(x)
$$
has $\Aut C = \{1,\iota\}$ if and only if
either $g\not\equiv 2\bmod 4$, or $g\equiv 2\bmod 4$ and
\begin{enumerate}
\item[(i)] $g-2$ is a $2$-power and $p(x+1)+p(x)\not\in\cH$;
\item[(ii)] $g-2$ is not a $2$-power and
$$p(x+1)+p(x)\not\in
      \cH + (x^4+x^2+1)((x+1)^{2g-3}+x^{2g-3})\neq \cH.$$
\end{enumerate}
\end{lemma}
\begin{proof}
Suppose $x\mapsto ax+b$  gives rise to a
non-extra automorphism $\lambda$ of $C$.

First we suppose $g\geq 7$.
If $b=0$ then (\ref{C5}) implies that $a=1$ and so $\lambda$ is not extra.
Otherwise, since $\deg(f(ax+b)+f(x))=2g$, all odd-degree terms in
$f(ax+b)+f(x)$ of degree $> g$ vanish.
Because $2g-5>g$ by our assumption,
the coefficients of $x^{2g-1}$, $x^{2g-3}$ and $x^{2g-5}$ are zero.
That is,
\begin{eqnarray}
\binom{2g+1}{2}b^2+1+a^2 = 0\\
\binom{2g+1}{4}b^4+\binom{2g-1}{2}b^2+1+a^4 = 0\\
\binom{2g+1}{6}b^4+\binom{2g-1}{4}b^2+\binom{2g-3}{2} = 0
\end{eqnarray}
Simplifying, we get respectively
\begin{eqnarray}
\label{E1}
gb^2+1+a^2=0\\
\label{E2}
\frac{g(g-1)}{2}b^4+(g-1)b^2+1+a^4 = 0\\
\label{E3}
\frac{g(g-1)(g-2)}{2}b^4+\frac{(g-1)(g-2)}{2}b^2+g =0.
\end{eqnarray}
Substitute (\ref{E1}) to (\ref{E2}) we get
$$\frac{g(g-1)}{2}b^4 + (g-1) b^2 + g^2 b^4 = 0;$$
and so
$$\frac{g(3g-1)}{2}b^2 + (g-1) = 0.$$ Thus $\frac{g(3g-1)}{2}=g-1$ and
$g\equiv 1,2\bmod 4$. But (\ref{E3}) implies $g\not\equiv 1\bmod 4$.

From now on we assume $g\equiv 2\bmod 4$.  Under this condition we get
$a=b=1$ by (\ref{E1}) and (\ref{E2}).  Once again, we use (\ref{C5})
to get
$$f(x+1)+f(x)=(p(x+1)+p(x))+ \gamma(x)   \in \cH,$$
where $\gamma(x) = (x^4+x^2+1)((x+1)^{2g-3}+x^{2g-3})$.

We claim that $\gamma(x)\in \cH$ if and only if $g-2$ is a $2$-power.
Suppose $\gamma(x)\in\cH$. We have $\deg(\gamma)=2g$ and
the its odd-degree terms are
\begin{eqnarray*}
&&\binom{2g-3}{2}x^{2g-1}\\
&+&(\binom{2g-3}{2}+\binom{2g-3}{4})x^{2g-3}\\
&+&(\binom{2g-3}{2}+\binom{2g-3}{4}+\binom{2g-3}{6})x^{2g-5}\\
&+&(\binom{2g-3}{4}+\binom{2g-3}{6}+\binom{2g-3}{8})x^{2g-7}\\
&+&\cdots\\
&+&(\binom{2g-3}{2m-4}+\binom{2g-3}{2m-2}+\binom{2g-3}{2m})x^{2g-(2m-1)}
\end{eqnarray*}
Set the odd-degree terms of degree $>g$ zero, and
use the identity
$\binom{2g-3}{2n}=\binom{g-2}{n}$ over $\F_2$ for all $n$, we have
\begin{eqnarray*}
\binom{g-2}{1}&=&0\\
\binom{g-2}{1}+\binom{g-2}{2}&=&0\\
\binom{g-2}{1}+\binom{g-2}{2}+\binom{g-2}{3}&=&0  \\
\binom{g-2}{2}+\binom{g-2}{3}+\binom{g-2}{4}&=&0  \\
\vdots\\
\binom{g-2}{m-2}+\binom{g-2}{m-1}+\binom{g-2}{m}&=&0
\end{eqnarray*}
for $m<\frac{g+1}{2}$. But we already have
$\binom{g-2}{1}=\binom{g-2}{2}=\binom{g-2}{3}=0$, so this system of
equations has a solution if and only if
$\binom{g-2}{m}=0$ for all $m\leq\frac{g}{2}$.
That is, $g-2$ is a $2$-power.
This proved parts (i) and (ii).

When $g=4$, we follow the same argument but only simpler. Namely,
any non-trivial automorphism $\lambda$ will lead to
(\ref{E1}) and (\ref{E2}) and hence $g\equiv 1,2\bmod
4$.  Contradiction.
\end{proof}

\begin{case}
Suppose $r=0$ and $g=4$ or $g\geq 7$:
\end{case}

Let $f(x) = x^{2g+1}+x^{2g-1}+x^{2g-3}+p(x)$, where $p(x)$ is any
polynomial in $\F_2[x]$ of degree $\leq 2g-6$  such that
$g\not\equiv 2\bmod 4$, or $g\equiv 2\bmod 4$ and
\begin{enumerate}
\item[(i)] if $g-2$ is a $2$-power, then
let $p=x^n+x^{n-2}+\mbox{(lower-degree terms)}$ where $n\equiv 3\bmod 4$; or
\item[(ii)] if $g-2$ is not a $2$-power, then let
$p\in \cH$.
\end{enumerate}

We shall verify our construction above.  If $g\not\equiv 2\bmod 4$ it
follows from Lemma~\ref{L4}. Suppose $g\equiv 2\bmod 4$. It can be
easily checked that part (i) implies $p(x+1)+p(x)\not\in\cH$ so it
follows from part (i) of the same Lemma.  In part (ii) $p\in\cH$ implies that
$p(x+1)+p(x)\in\cH$.  Since $g-2$ is not a $2$-power, $\cH
+(x^4+x^2+1)((x+1)^{2g-3}+x^{2g-3})$ is a non-trivial coset of $\cH$,
hence is disjoint from $\cH$. So part (ii) follows from part (ii) of
the same Lemma again.

\begin{case}
Suppose $r=0$ and $g=6$:
\end{case}

Let $f = x^{2g+1} + x^{2g-3} + x^{2g-5}+p(x)$ where $p(x)$ is a
polynomial in $\F_2[x]$ of degree $\leq 2g-6$.

Suppose $x\mapsto ax+b$ gives rise to an automorphism $\lambda$ of $C$.
For any $g\equiv 2\bmod 4$ we show that the only possible extra
automorphism is the one given by $a,b$ are both $3$-rd roots of unity
over $\F_2$. Apply~(\ref{C5}) to coefficients of
$x^{2g}, x^{2g-1}, x^{2g-3}, x^{2g-5}$, those are
$a^{2g}b, 1+a^{2g-3}(1+b^4), 1+a^{2g-5}$.
If $b=0$ then $a=1$ so it is trivial.
If $b\neq 0$ then $a=b+1$ and $a^3=1$. If it is not trivial then
$a,b$ are $3$-rd roots of unity.

When $g=6$ we have $a^{2g-5}=a^7=1$ and $a^6=1$ so $a=1$.
This implies $b=0$.  So $\lambda$ is trivial.

\begin{case}
Suppose $r=0$ and $g=3$ or $5$:
\end{case}

Let $f = x^{2g+1}+x^{2g-3}+p(x)$ where $p(x)$ is a polynomial
in $\F_2[x]$ of degree $2g-5$. In fact, this construction works
for every odd $g\geq 3$.

Suppose $x\mapsto ax+b$ gives rise to an automorphism $\lambda$ of $C$.
The coefficient of $x^{2g}$ and $x^{2g-1}$ in $f(ax+b)+f(x)$
are $a^{2g}b$ and $a^{2g-1}b$, respectively.
At least one of them has to vanish by (\ref{C5}), so $b=0$.
This implies $a=1$ by applying (\ref{C5}) again.

\begin{acks}
The author thanks R. Coleman and B. Poonen for valuable
correspondences and remarks on draft of this paper. She is
grateful to H. Lenstra for support and the Mathematisch Institute
of Universiteit Leiden for its hospitality, where part of this
work was done in 2000. Finally she thanks the referee for
comments.
\end{acks}

\bibliographystyle{amsplain}

\end{document}